\documentclass[11pt]{article}
\usepackage{amsfonts}
\usepackage{amssymb}
\usepackage{graphicx}
\usepackage{amsmath}
\renewcommand{\raggedright}{\leftskip=0pt \rightskip=0pt plus 0cm}
\raggedright

\renewcommand{\raggedright}{\leftskip=0pt \rightskip=0pt plus 0cm}
\raggedright

\def\N{{\mathbb N}}

\def\R{{\mathbb R}}

\def\bb{\begin}
\def\bc{\begin{center}}
\def\ec{\end{center}}
\def\be{\begin{equation}}
\def\ee{\end{equation}}
\def\ba{\begin{array}}
\def\ea{\end{array}}
\def\bea{\begin{eqnarray}}
\def\eea{\end{eqnarray}}
\def\beaa{\begin{eqnarray*}}
\def\eeaa{\end{eqnarray*}}

\def\hh{\!\!\!\!}
\def\EM{\hh &   &\hh}
\def\EQ{\hh & = & \hh}
\def\EE{\hh & \equiv & \hh}

\def\GE{\hh & \ge & \hh}

\def\AND#1{\hh & #1 & \hh}
\def\DEQ{\AND{:=}}

\def\al{\alpha}
\def\bt{\beta}
\def\la{\lambda}

\def\vp{\varphi}

\def\da{\delta}

\def\ck{\check}

\def\d{\cdot}

\def\dd{\cdots}
\def\oo{\infty}
\def\f{\frac}
\def\pa{\partial}
\def\z{\left}
\def\y{\right}
\def\q{\quad}
\def\qq{\qquad}

\def\tm{\times}

\def\rd{\,{\rm d}}

\def\du{\,{\rm d}u}
\def\dv{\,{\rm d}v}
\def\dx{\,{\rm d}x}
\def\dy{\,{\rm d}y}

\def\E{{\mathcal E}}
\def\mcl{{\mathcal L}}

\def\Lo{\mcl^1}

\def\Li{\mcl^\infty}

\def\nlo{\|\cdot\|_1}

\def\nli{\|\cdot\|_\infty}

\def\nmt#1{\|#1\|_2}

\def\ii{\int_I}

\def\lto{\rightarrow}

\def\to{\lto}

\def\andq{\quad \mbox{ and } \quad}
\def\qqf{\qquad \forall}

\def\Fr{Fr\'echet\ }

\def\ifl{\iffalse}
\def\Proof{\noindent{\bf Proof} \quad}
\def\qed{\hfill $\Box$ \smallskip}

\def\nn{\nonumber}
\def\lb{\label}
\def\x#1{{\rm (\ref{#1})}}
\def\wx#1#2#3#4#5#6#7#8{\bibitem{#1} {#2}, {#3}, {\it #4}, {{\bf #5}} {(#6)}, {#7}--{#8}.}

\def\inn#1#2{\z\langle #1, \, #2 \y\rangle}

\def\wrt{with respect to\ }

\begin{document}

\title{On the Second-order Fr\'echet Derivatives of Eigenvalues of Sturm-Liouville Problems in Potentials}

\author{Shuyuan Guo, \qquad Guixin Xu, \qquad Meirong Zhang\footnote{Correspondence author. This author is supported by the National Natural Science Foundation of China (Grant no. 11790273).}\\
\bigskip
Department of Mathematical Sciences, Tsinghua University, Beijing 100084, China\\
\bigskip
E-mails: {\tt guo-sy17@mails.tsinghua.edu.cn} (S. Guo)\\
{\tt guixinxu$_-$ds@163.com} (G. Xu)\\
{\tt zhangmr@tsinghua.edu.cn} (M. Zhang) }

\date{}%

\maketitle

\begin{abstract}
The works of V. A. Vinokurov have shown that eigenvalues and normalized eigenfunctions of Sturm-Liouville problems are analytic in potentials, considered as mappings from the Lebesgue space to the space of real numbers and the Banach space of continuous functions respectively. Moreover, the first-order Fr\'echet derivatives are known and paly an important role in many problems. In this paper, we will find the second-order Fr\'echet derivatives of eigenvalues in potentials, which are also proved to be negative definite quadratic forms for some cases. 
\end{abstract}

\bigskip

{\small {\bf Mathematics Subject Classification (2010)}:
34B24; 34L10; 49J50; 58C20.

}

\section{Introduction}
\setcounter{equation}{0} \lb{first}

It is an important issue and has a long history to study the dependence of eigenvalues of differential operators on coefficients, boundary data and domains involved in the problems \cite{CH53, KWZ97, KWZ99, MZ96}. In this paper, we are concerning with the dependence of eigenvalues on potentials. More precisely, let $q=q(x)\in \Lo:=L^1(I,\R)$ be an integrable potential, where $I=[0,\ell]$, $\ell>0$ is a closed interval. We consider the Sturm-Liouville eigenvalue problem
    \be \lb{Eq}
    -z'' + q(x) z=\la z, \qq x\in I,
    \ee
with the boundary conditions
    \bea\lb{bc0}
    z(0) \cos \al -z'(0)\sin \al \EQ 0, \\
    \lb{bc1}
    z(\ell) \cos \bt -z'(\ell)\sin \bt \EQ 0,
    \eea
where $\al\in[0,\pi)$ and $\bt \in(0,\pi]$ are given parameters. It is well-known that problem \x{Eq}---\x{bc1} has a discrete spectrum consisting of an increasing infinite sequence of (real, simple) eigenvalues $\la_n$ such that $\la_n\to+\oo$ as $n\to \oo$. See, for example, \cite{Ze05}. The eigenvalues $\la_n$ depend on potential $q$ and on boundary data $(\al,\bt)$ as well.

Let us fix $(\al,\bt)$ and consider $\la_n=\la_n(q)$ as (nonlinear) functionals of potentials $q\in \Lo$. Associated with $\la_n$ is the corresponding eigenfunction $E_n(x)=E_n(x;q)$, normalized as
    \be \lb{En1}
    \nmt{E_n} = \z( \ii E_n^2(x)\dx\y)^{1/2}=1.
    \ee
Moreover, in order that $E_n(x)$ is uniquely determined, it is required that
    \be \lb{En2}
    E_n(x)>0 \mbox{ on some neighborhood of the type }(0,\da).
    \ee
Among many studies on Sturm-Liouville problems, around 2005, Vinokurov and his collaborator \cite{V05b, V05c, VS05} have systematically studied the dependence of solutions and eigenvalues on potentials. One of their results in \cite{V05c} is that when $\Lo$ and $C(I,\R)$ are respectively endowed with the $L^1$ norm $\|\d\|_1$ and the supremum norm $\|\d\|_\oo$, the nonlinear mappings
    \be \lb{mps}
    \z\{ \ba{lcl} (\Lo,\nlo) \to \R,& & q\mapsto \la_n(q),\\
    (\Lo,\nlo) \to (C(I,\R),\|\d\|_\oo), & & q\mapsto E_n(\d;q),
    \ea \y.
    \ee
are proved to be analytic in the sense of \cite{B67, B71}. Moreover, the \Fr derivative of eigenvalue $\la_n(q)$ is given by
    \be \lb{D11}
    \pa_q \la_n(q)(h):=\z.\pa_s \la_n(q+s h)\y|_{s=0} \equiv \ii \z(E_n(x;q)\y)^2 h(x) \dx
    \ee
for $h \in \Lo$. The \Fr derivative of eigenfunction $E_n(\d;q)$ in $q\in \Lo$ can also be found from \cite{V05c}. In fact, result \x{D11} for the Dirichlet eigenvalues and potentials in $L^2(I,\R)$ was already obtained in \cite{PT87}. As for the first Dirichlet eigenvalue of the Laplacian with potentials, similar formula of the \Fr derivative is derived very recently in \cite{IV19}. Formula \x{D11} can also be written as
    \be \lb{D12}
    \pa_q \la_n(q)=\z(E_n(\d;q)\y)^2,
    \ee
understood as a kernel function in the dual space $(\Lo,\nlo)^*=L^\oo(I,\R)$.

One simple implication of \x{D12} is that eigenvalues $\la_n(q)$ are strictly increasing in $q\in \Lo$, because the derivatives are positive. As for the continuous dependence, about ten years ago, one of the authors of this paper and his collaborator have obtained a further result. That is, when the norm topology $\nlo$ in $\Lo$ is replaced by the topology $w_1$ of weak convergence in $\Lo$, the mappings in \x{mps} are still continuous \cite{MZ10, Zh08}. Based on such a complete continuity of eigenvalues in potentials, the derivatives \x{D12} of eigenvalues are applied in \cite{WMZ09, Zh09} to study some typical optimization problems on eigenvalues of Sturm-Liouville operators. This leads to some connection between the linear and nonlinear stationary Schr\"odinger equations. It is interesting to note that in very recent papers \cite{IV18, IV19}, such a connection is also established for some inverse spectral problems.

Because of these applications, it is convincing that the higher-order \Fr derivatives of eigenvalues and eigenfunctions will be useful. However, though $\la_n(q)$ and $E_n(x;q)$ are analytic in $q\in\Lo$, as far as we know, even their second-order \Fr derivatives in potentials are not available in the literature. The aim of this paper is to meet this gap.

For $q, \ h\in \Lo$, let us define the second-order \Fr derivative \cite{Ze85} of $\la_n(q)$ by
    \[
\pa_{qq} \la_n(q)(h):=\z.\pa_{ss} \la_n(q+s h)\y|_{s=0}\in \R.
    \]
For simplicity, denote
    \be\lb{L}
    L:=\pa_q \la_n(q)(h)=\ii \z(E_n(x;q)\y)^2 h(x) \dx\in \R.
    \ee
We will obtain the following results.

    \bb{thm} \lb{A}
For $q, \ h\in \Lo$, let $L$ be as in \x{L} and $U_n(x)= U_n(x;q,h)$ be the unique solution of the following inhomogeneous linear ODE
    $$
    -z'' +\z(q(x)-\la_n(q)\y) z = -E_n(x)\z(h(x)-L\y),\q x\in I, \eqno(E)
    $$
satisfying the initial conditions
    $$
    z(0)=z'(0)=0.\eqno (I)
    $$
Then the second-order \Fr derivative of eigenvalue is given by
    \be \lb{D21}
\pa_{qq} \la_n(q)(h)=2\ii E_n(x)\z(h(x)-L\y) U_n(x) \dx.
    \ee
Moreover, it can also be expressed as a quadratic form of $h$
    \be \lb{D22}
    \pa_{qq} \la_n(q)(h)=\int_{I^2} J_n(x,y) h(x) h(y) \dx\dy.
    \ee
Here $U_n(x;q,h)$ and $J_n(x,y)=J_n(x,y;q)$ are explicitly given in formula \x{Un3} and in \x{Jn1}-\x{Jn2} respectively.
    \end{thm}

Next, we are still using the solutions $U_n(x)$, but with some restriction on the boundary condition \x{bc1} at the right end-point $x=\ell$.

    \bb{thm} \lb{B}
For any $q, \ h\in \Lo$, let $U_n(x)$ be as in Theorem \ref{A}. Then

{\rm (i)} $U_n(x)$ also satisfies the boundary conditions \x{bc0}-\x{bc1}.

{\rm (ii)} Assume that \x{bc1} takes the following boundary condition
    \be \lb{bc2}
    z(\ell)=0 \qq \mbox{or} \qq z'(\ell)=0.
    \ee
Then, in this case, the second-order \Fr derivative can also be expressed as
    \be \lb{D23}
\pa_{qq} \la_n(q)(h)=-2\ii \z(U_n'^2 +\z(q(x)-\la_n(q)\y) U_n^2\y)\dx.
    \ee
    \end{thm}

The paper is organized as follows. The proofs of Theorem \ref{A} and Theorem \ref{B} will be given in \S \ref{second} and in \S \ref{conc} respectively. As for the first Dirichlet eigenvalues $\la_1^D(q)$ of \x{Eq}, we will use \x{D23} to prove in Theorem \ref{C} that $\pa_{qq} \la_1^D(q)(h)$ is a negative definite quadratic form of $h\in \Lo$. This result is consistent with the concavity of $\la_1^D(q)$ in $q\in \Lo$. See the discussion in Remark \ref{con}. At the end of the paper, we will propose some further problems on the second-order \Fr derivatives.

\section{The Second-order \Fr Derivatives of Eigenvalues}
\setcounter{equation}{0} \lb{second}

Let $q\in \Lo$. For $h\in \Lo$ and $s\in \R$, let us denote
    \be \lb{Qs}
    Q(x,s):=q(x)+ s h(x)-\la_n(q+s h)\in \Lo.
    \ee
We consider the linear ODE
    \be \lb{Eq1}
    -z'' +Q(x,s) z=0, \qq x\in I.
    \ee
Here $'=\f{\rd}{\dx}$ and $'' =\f{\rd^2}{\dx^2}$ are also written as $\pa_x$ and $\pa_{xx}$ respectively. For each $s$, let $z=\vp(x,s)$ be the solution of Eq. \x{Eq1} satisfying the initial condition
    \be \lb{i0}
(z(0),z'(0))=(\sin \al, \cos \al).
    \ee
That is, for each $s$, one has
    \bea \lb{Eq2}
    \EM -\pa_{xx}\vp(x,s) + Q(x,s)\vp(x,s) = 0,\\
    \lb{I0}
    \EM (\vp(0,s),\pa_{x}\vp(0,s)) = (\sin \al, \cos \al).
    \eea
Moreover, from the definition of eigenvalues, for each $s$, $\vp(\d,s)$ also satisfies the boundary condition \x{bc1}
    \be \lb{Bl1}
    \vp(\ell,s)\cos \bt -\pa_x\vp(\ell,s)\sin \bt = 0.
    \ee

Since $\la_n(q+s h)$ is analytic in $s$, as a mapping from $\R$ to $\Lo$, $Q(\d,s)$ is then analytic in $s$. Moreover, since the solution of initial value problem of linear ODE is analytic in coefficient potential \cite{V05b}, $\vp(x,s)$ and $\pa_x \vp(x,s)$ are also analytic in $s$. Thus, by differentiating Eq. \x{Eq2} \wrt $s$ twice, we obtain the following inhomogeneous linear ODEs
    \bea \lb{E31}
    \EM -\pa_{xx}\pa_{s}\vp(x,s) + Q(x,s)\pa_{s}\vp(x,s) = -\pa_{s} Q(x,s)\d\vp(x,s),\\
    \lb{E32}
    \EM -\pa_{xx}\pa_{ss}\vp(x,s) + Q(x,s)\pa_{ss}\vp(x,s)\nn\\
    \EQ -\pa_{ss} Q(x,s)\d\vp(x,s) - 2\pa_{s} Q(x,s)\d\pa_{s}\vp(x,s).
    \eea
Let us consider these equations at $s=0$. For simplicity, denote
    \be\lb{L-M}
    \z\{ \ba{l} L := \z.\pa_s\la_n(q+s h)\y|_{s=0}=\pa_q \la_n(q)(h)\in \R,\\
    M := \z.\pa_{ss}\la_n(q+s h)\y|_{s=0}=\pa_{qq} \la_n(q)(h)\in \R,\\
    Q(x):= Q(x,0) = q(x) -\la_n(q),\\
    \vp(x):= \vp(x,0),\\
    \vp_k(x) := \z. \pa_{s^k} \vp(x,s)\y|_{s=0},\qq k=1,2.
    \ea \y.
    \ee
By \x{Qs}, one has
    \(
    \z. \pa_s Q(x,s)\y|_{s=0}=h(x)-L
    \)
and
    \(
    \z. \pa_{ss} Q(x,s)\y|_{s=0}=- M.
    \)
Thus Eq. \x{E31} and Eq. \x{E32} mean that $z=\vp_k(x)$ is a solution of the inhomogeneous linear ODE
    \be \lb{Eq3}
    -z'' + Q(x) z = f_k(x), \qq x\in I,
    \ee
where
    \be\lb{f1x}
    f_1(x)= (L-h(x))\vp(x) \andq
    f_2(x)= M \vp(x)-2(h(x) - L)\vp_1(x).
    \ee
Moreover, by differentiating initial condition \x{I0} and boundary condition \x{Bl1} \wrt $s$ twice, we know that $z=\vp_k(x)$, $k=1,2$ satisfy
    \bea\lb{I1}
    \EM \z(z(0),z'(0)\y) = (0,0),\\
    \lb{bc11}
    \EM z(\ell)\cos \bt -z'(\ell)\sin \bt = 0.
    \eea

    \bb{lem}\lb{S}
Let $q\in \Lo$ be given. For any $h\in \Lo$, one has

{\rm (i)} $\vp_k(x)=\vp_k(x;q,h)$, $k=1,2$ are uniquely determined. Actually, they are the solutions of  the initial value problems \x{Eq3}-\x{I1}.

{\rm (ii)} $\vp_k(x)=\pa_{s^k} \vp(x,s)|_{s=0}$, $k=1,2$ are given by
    \be \lb{vpk}
    \vp_k(x)= \int_0^x W(x,y)f_k(y)\dy , \qq x\in I,
    \ee
where $f_1$ and $f_2$ are as in \x{f1x}, and $W(x,y)$ is as in \x{W} below. In particular,
    \be \lb{vp1}
    \vp_1(x)= \int_0^x W(x,y)\vp(y)(L-h(y))\dy, \qq x\in I.
    \ee

{\rm (iii)} $z=\vp(x), \ \vp_1(x), \ \vp_2(x)$ satisfy the boundary conditions \x{bc0}-\x{bc1}.
    \end{lem}

\Proof (i) This is clear from the above deductions.

(ii) Let $\psi_i(x)=\psi_i(x;q)$, $i=1,2$ be the fundamental solutions of the homogeneous linear ODE
    \be \lb{Z0}
    -z''+ Q(x) z= -z''+ (q(x)-\la_n(q)) z=0,
    \ee
i.e. the solutions of Eq. \x{Z0} satisfying the initial conditions $(z(0),z'(0))=e_1:=(1,0)$ and $(z(0),z'(0))=e_2:=(0,1)$ respectively. Define
    \be\lb{W}
    W(x,y):= \z| \ba{cc} \psi_1(x) & \psi_2(x) \\ \psi_1(y) & \psi_2(y)\ea \y|=\psi_1(x)\psi_2(y)- \psi_2(x)\psi_1(y).
    \ee
By applying the formula-of-constant-variant to \x{Eq3}-\x{I1}, we know that $\vp_k(x)$ is given by \x{vpk}. In particular, \x{vp1} follows from \x{f1x} and \x{vpk}. Here one can notice that $\psi_i(x)$ and $W(x,y)$ depend only on $q$, not on $h$.

(iii) Due to the choice \x{i0} of the initial conditions for $\vp(x)$, $\vp(x)$ satisfies boundary condition \x{bc0}. Moreover, by the definition of eigenvalues $\la_n(q)$, one knows that $\vp(x)$ satisfies \x{bc1} as well.

For $k=1,2$, it is clear from \x{I1} and \x{bc11} that $z=\vp_k(x)$ must satisfy the boundary conditions \x{bc0}-\x{bc1}. \qed

Lemma \ref{S} (iii) shows that $\vp(x)$ is an eigenfunction associated with $\la_n(q)$. One then sees that the normalized eigenfunction satisfying \x{En1} and \x{En2} is
    \be \lb{En3}
    E_n(x)=E_n(x;q) \equiv {\vp(x)}/{\nmt{\vp}}.
    \ee

To derive the formulas for the \Fr derivatives $L$ and $M$, we can exploit the Fredholm principle \cite{DS58}.

    \bb{lem} \lb{FP}
Consider the following inhomogeneous linear ODE
    \be \lb{Z2}
    -z'' +Q(x) z = f(x),\qq x\in I,
    \ee
where $f(x)\in \Lo$. If Eq. \x{Z2} admits a solution $z(x)$ satisfying the boundary conditions \x{bc0}-\x{bc1}, it is necessary that
    \be \lb{Z20}
    \ii \vp(x) f(x)\dx =0.
    \ee
    \end{lem}

\Proof Let $\vp(x)$ be as in \x{L-M}. From Lemma \ref{S}  (iii), $z=\vp(x)$ is a solution of boundary value problem \x{Z0}-\x{bc0}-\x{bc1}, where Eq. \x{Z0} is the corresponding homogeneous linear ODE of Eq. \x{Z2}. Hence the solvability of problem \x{Z2}-\x{bc0}-\x{bc1} is actually equivalent to condition \x{Z20}. \qed

{\bf Proof of  Theorem \ref{A}} Applying Lemma \ref{FP} to Eq. \x{Eq3}, we obtain
    \be \lb{Z22}
    \ii \vp(x) f_k(x)\dx =0,\qq k=1,2.
    \ee
For $k=1$, it follows from \x{f1x} that equality \x{Z22} is
    \(
    \ii \z(L-h(x)\y)\vp^2(x)\dx =0,
    \)
i.e.
    \be \lb{Lh}
    L = \f{\ii \vp^2(x) h(x)\dx}{\ii \vp^2(x) \dx}= \ii \z(\f{\vp(x)}{\nmt{\vp}} \y)^2 h(x)\dx= \ii E_n^2(x) h(x)\dx.
    \ee
See \x{En3}. This gives another deduction of formula \x{D11}, which is different from that in \cite{V05c}.

Let now $L$ be as in \x{Lh}. For $k=2$, it follows from \x{f1x} that equality \x{Z22} is
    \[
    \ii \z(M \vp^2(x)-2\vp(x)(h(x) - L)\vp_1(x)\y)\dx =0.
    \]
By using \x{En3}, this yields
    \[
    M =2\ii E_n(x)(h(x) -L) \f{\vp_1(x)}{\nmt{\vp}}\dx.
    \]
Thus we have obtained \x{D21}, where
    \be \lb{Un1}
    U_n(x)= U_n(x;q,h):=\f{\vp_1(x)}{\nmt{\vp}} \equiv \f{\z.\pa_s \vp(x,s)\y|_{s=0}}{\nmt{\vp}}.
    \ee
Dividing Eq. \x{Eq3} by the factor $\nmt{\vp}$ and making use of conditions \x{f1x}, \x{I1}, one sees that $U_n(x)$  is just the solution of the initial value problem $(E)$-$(I)$. From \x{vp1}, \x{En3} and \x{Un1}, it is easy to see that
    \be \lb{Un3}
    U_n(x)=U_n(x;q,h)\equiv \int_0^x W(x,y)E_n(y)(L-h(y))\dy, \qq x\in I,
    \ee
where $W(x,y)$ and $L$ are in \x{W} and \x{Lh} respectively.

Now we are using formulas \x{D21} and \x{Un3} to deduce formula \x{D22}. Define
    \bea \lb{hvp0}
    \hat h(x) \DEQ E_n(x) (h(x)-L)\\
    \lb{hvp}
    \EQ \ii E_n(x) E_n^2(u)(h(x)- h(u)) \du,
    \eea
because
    \(
    \ii E_n^2(u) \du =1.
    \)
By \x{Un3}, one has
    \[
    U_n(x)=-\int_0^x W(x,y) \hat h(y) \dy.
    \]
Hence \x{D21} gives
    \bea \lb{D4}
    M \EQ 2\int_0^\ell \hat h(x) U_n(x) \dx=2 \int_0^\ell \hat h(x) \z(-\int_0^x W(x,y) \hat h(y) \dy \y)\dx\nn\\
    \EQ \int_{I^2} G(x,y) \hat h(x) \hat h(y) \dx\dy.
    \eea
Here $G(x,y): I^2 \to\R$ is the following symmetrization of $W(x,y)$
    \be \lb{G}
    G(x,y) :=\left\{
    \begin{aligned}
    W(x,y) & \;\;\;0\leq x\le y\le \ell, \\
    -W(x,y) & \;\;\;0\le y\le x\le \ell.
    \end{aligned}
    \right.
    \ee
Obviously, we have $G(x,y)=G(y,x)$, i.e. $G(x,y)$ is symmetric.

By \x{hvp0}, we have
    \beaa
    \hat h(x)\hat h(y)\EQ E_n(x) E_n(y) h(x)h(y) -E_n(x) E_n(y)\d L h(x)\\
    \EM -E_n(x) E_n(y)\d L h(y) + E_n(x) E_n(y)\d L\d L.
    \eeaa
Denote
    \[
    \E(x,y,u,v):= E_n(x) E_n(y) E_n(u) E_n(v).
    \]
Then \x{D4} is
    \[
    M=M_1-M_2-M_3+M_4,
    \]
where, by using \x{hvp},
    \beaa
    M_1 \EQ \int_{I^2} \z(\int_{I^2}G(x,y) E_n^2(u)E_n^2(v)\du\dv\y) E_n(x) E_n(y) h(x) h(y)\dx\dy\\
    \EQ \int_{I^2} \z(\int_{I^2}G(x,y) E_n(u)E_n(v) \E(x,y,u,v)\du\dv\y) h(x) h(y)\dx\dy,\\
    M_2 \EQ \int_{I^2} G(x,v) E_n(x) E_n(v) \d L h(x)\dx\dv\\
    \EQ \int_{I^2} G(x,v) E_n(x) E_n(v) \z(\ii E_n^2(y) h(y) \dy\y) h(x)\dx\dv\\
    \EQ \int_{I^2} \z(\int_{I^2} G(x,v) E_n(u) E_n(y) \E(x,y,u,v)\du \dv\y)h(x) h(y)\dx\dy.
    \eeaa
Similarly,
    \beaa
    M_3 \EQ \int_{I^2} \z(\int_{I^2} G(u,y) E_n(v)E_n(x) \E(x,y,u,v)\du\dv\y) h(x) h(y)\dx\dy,\\
    M_4 \EQ \int_{I^2} \z(\int_{I^2} G(u,v) E_n(x) E_n(y) \E(x,y,u,v)\du\dv\y)h(x) h(y)\dx\dy.
    \eeaa
Thus, by defining $J_n(x,y,u,v): I^4 \to \R$ as
    \bea\lb{Jn1}
    J_n(x,y,u,v)\DEQ G(x,y) E_n(u)E_n(v) - G(x,v) E_n(u) E_n(y)\nn\\
    \AND{-}G(u,y) E_n(x) E_n(v)+G(u,v) E_n(x) E_n(y),
    \eea
and by defining $J_n(x,y): I^2 \to \R$ as
    \be\lb{Jn2}
    J_n(x,y) :=\z(\int_{I^2}J_n(x,y,u,v)E_n(u)E_n(v)\du\dv \y) E_n(x) E_n(y),
    \ee
we know that $M$ is expressed as the integral form \x{D22}.

The proof of Theorem \ref{A} is completed.    \qed

    \bb{rem}\lb{ker} {
The kernel $J_n(x,y,u,v)$ of \x{Jn2}, a continuous function defined on $I^4$, is determined from \x{G} and \x{Jn1}. These kernels have the following symmetries
    \[
    J_n(x,y,u,v) \equiv J_n(u,v,x,y)\andq J_n(x,y) \equiv J_n(y,x).
    \]
\ifl Moreover,
    \[
    \int_{I^4} J_n(x,y,u,v)\du\dv \dx\dy =0.
    \]
In particular, $J_n(x,y,u,v)$ must be sign-changing on $I^4$. \fi
}
    \end{rem}

\ifl

\bb{exa} \lb{q0} {\rm Let $\ell=\pi$ and $q(x)\equiv 0$ in Eq. \x{Eq}. We consider the Dirichlet eigenvalues $\la_n(0)=n^2$, $n\in \N$. Then
    \[
    \ba{ll} \psi_1(x)=\cos (n x), & \psi_2(x)=\f{\sin (n x)}{n},\\
    G(x,y) = \f{\sin (n|x-y|)}{n},& E_n(x)= \sqrt{\f{2}{\pi}}\sin (n x).\ea
    \]
From these, one can obtain, for $0\le y\le x\le \pi$,
    \beaa
    J_n(x,y)\EQ J_n(y,x)\\
    \EQ G(x,y)E_n(x)E_n(y)-\z(\int_{[0,\pi]} G(x,v) E_n(v)\dv\y) E_n(x)E_n^2(y)\\
    \EM - \z(\int_{[0,\pi]} G(u,y) E_n(u)\du\y) E_n^2(x)E_n(y)\\
    \EM + \z(\int_{[0,\pi]^2} G(u,v) E_n(u)E_n(v)\du\dv\y) E_n^2(x)E_n^2(y)\\
    \EE \f{2\sin(n x) \sin(n y)}{(n\pi)^2}\Bigl(2n x\cos(n x) \sin(n y) +2 n y\cos(n y) \sin(n x)\\
    \EM -\z(2n\pi \cos(n x) +\sin(n x) \y)\sin(n y) \Bigr).
    \eeaa
This kernel changes sign on the domain $[0,\pi]^2$.\qed}
\end{exa}

\fi

\section{The Concavity of Eigenvalues in Potentials}
\setcounter{equation}{0} \lb{conc}

We will derive formula \x{D23} for the second-order \Fr derivatives of eigenvalues in potentials.

{\bf Proof of Theorem \ref{B}} (i) For general boundary conditions \x{bc0}-\x{bc1}, it is clear from Lemma \ref{S} (iii) that the solution $U_n(x)$ of $(E)$ also satisfies the boundary conditions \x{bc0}-\x{bc1}.

(ii) Recall that $U_n(x)$ satisfies ODE
    \be \lb{Une}
    -U''_n + Q(x) U_n =-E_n(x) (h(x)-L).
    \ee
Since we are now considering the boundary conditions \x{bc0}-\x{bc2}, from the proof of Lemma \ref{S}, we know that $U_n(x)=\vp_1(x)/\nmt{\vp}$ satisfies
    \be \lb{bc3}
    U_n(0)=U'_n(0)=0\andq U_n(\ell) U'_n(\ell)=0.
    \ee

Multiplying Eq. \x{Une} by $U_n(x)$ and then integrating on $I$, we obtain
    \beaa
    \EM \ii E_n(x) (h(x)-L) U_n(x) \dx
    = -\ii\z(-U_n''+Q(x) U_n\y)U_n\dx \\
    \EQ \z. U_n(x)U'_n(x)\y|_{x=0}^\ell -\ii \z(U_n'^2 +Q(x) U_n^2\y)\dx= -\ii \z(U_n'^2 +Q(x) U_n^2\y)\dx,
    \eeaa
because we have \x{bc3}. Hence formula \x{D23} can be deduced from \x{D21} and the above equality. \qed

    \bb{rem} \lb{Vn}
{
In boundary conditions \x{bc0}-\x{bc1}, if \x{bc0} is restricted to be either $z(0)=0$ or $z'(0)=0$, dually we have
    \[
\pa_{qq} \la_n(q)(h)=-2\ii \z(V_n'^2 +\z(q(x)-\la_n(q)\y) V_n^2\y)\dx.
    \]
Here $z=V_n(x)$ is the solution of $(E)$ satisfying the initial conditions $z(\ell)=z'(\ell)=0$.}
    \end{rem}

As an example, let us consider the Dirichlet boundary conditions
    $$
    z(0)=z(1)=0. \eqno(D)
    $$
For $q\in\Lo$, we use $\la_n^D(q)$, $n\in \N$ to denote the eigenvalues of problem \x{Eq}-$(D)$. It is known from \cite{CH53} that $\la_1^D(q)$ has the following minimization characterization
    \be \lb{la1}
    \la_1^D(q) =\min_{z\in H_0^1(I), \ z\ne 0}\f{\ii (z'^2 + q(x) z^2) \dx}{\ii z^2 \dx}.
    \ee
For any $n\ge 2$, $\la_n^D(q)$ has the following maximin characterization
    \be \lb{lan}
    \la_n^D(q) =\max_{W_{n-1}}\min_{z\in W_{n-1}^\bot, \ z\ne 0}\f{\ii (z'^2 + q(x) z^2) \dx}{\ii z^2 \dx},
    \ee
where the maximum is taken over all subspaces $W_{n-1}$ of $H_0^1(I)$ of dimension $n-1$, and
    \[
    W_{n-1}^\bot:=\z\{z\in H_0^1(I): \ii w(x) z(x) \dx=0\mbox{ for all } w \in W_{n-1} \y\}.
    \]


    \bb{thm}\lb{C}
For the first Dirichlet eigenvalues $\la_1^D(q)$, the second-order \Fr derivatives satisfy $\pa_{qq} \la_1^D(q)(h)\le 0$ for all $h\in \Lo$.
    \end{thm}

\Proof Since the Dirichlet boundary conditions satisfy \x{bc0}-\x{bc2}, we can use formula \x{D23} for $\pa_{qq} \la_1^D(q)(h)$. By Theorem \ref{B}, $U_1(x)=U_1^D(x)$ satisfies $(D)$, i.e. $U_1(x)\in H_0^1(I)$. From the minimization characterization \x{la1}, one has
    $$
    \ii \z(U_1'^2 +q(x) U_1^2\y)\dx\ge \la_1^D(q)\ii U_1^2 \dx.
    $$
Thus \x{D23} shows that $\pa_{qq} \la_1^D(q)(h)\le 0$. \qed


    \bb{rem}\lb{con} {
{\rm (i)} It is standard from nonlinear analysis \cite{Ze85} that the negative definiteness of the second-order \Fr derivatives is the same as the concavity. As a result, Theorem \ref{C} implies that
    \be \lb{la11}
    \la_1^D\z(\tau q_1 +(1-\tau) q_2\y) \ge \tau \la_1^D(q_1) + (1-\tau) \la_1^D(q_2)
    \ee
for all $q_i\in \Lo$ and $\tau\in [0,1]$.

{\rm (ii)} The concavity \x{la11} of the first eigenvalue $\la_1^D(q)$ in $q$ can also be directly deduced from the minimization characterization \x{la1}. In fact, one has
    \bea \lb{la12}
    \EM \la_1^D\z(\tau q_1 +(1-\tau) q_2\y) \nn\\
    \EQ \min_{z\in H_0^1(I) \atop \nmt{z}=1}\ii \z(z'^2 + \z(\tau q_1(x) +(1-\tau) q_2(x)\y) z^2\y) \dx\nn\\
    \EQ \min_{z\in H_0^1(I) \atop \nmt{z}=1}\z( \tau \ii (z'^2 + q_1(x) z^2) \dx+ (1-\tau) \ii (z'^2 + q_2(x) z^2) \dx\y)\nn\\
    \GE \tau \d \min_{z\in H_0^1(I) \atop \nmt{z}=1} \ii (z'^2 + q_1(x) z^2) \dx
    +(1-\tau)\d\min_{z\in H_0^1(I) \atop \nmt{z}=1}  \ii (z'^2 + q_2(x) z^2) \dx \nn\\
    \EQ \tau \la_1^D(q_1) + (1-\tau) \la_1^D(q_2).
    \eea

{\rm (iii)} For the zeroth Neumann eigenvalue $\la_0^N(q)$ of problem \x{Eq}, it can be proved that $\pa_{qq} \la_0^N(q)(h)\le 0$ for all $h\in \Lo$.
}
    \end{rem}

We end the paper by two problems.

1. Arguing as in the deduction of \x{la12}, one can use the maximin characterization \x{lan} to obtain the concavity of $\la_n^D(q)$ in $q\in \Lo$. Therefore $\pa_{qq} \la_n^D(q)(h)$ is also negative definite. It is an interesting problem to give a direct proof for the negative definiteness of $\pa_{qq} \la_n^D(q)(h)$ for the case $n\ge 2$.

2. We have known from \cite{V05c} that eigenfunctions $E_n(x;q)$ are also analytic in $q\in \Lo$. It is then an important problem to find the second-order \Fr derivatives of $E_n(x;q)$ in $q\in \Lo$.

\section*{Acknowledgments}

The authors would like to thank Professors Irina V. Astashova and Alexey V. Filinovskiy for their interests and suggestions in this work.

\end{document}